\documentclass[10pt]{amsart}

\usepackage[latin1]{inputenc}
\usepackage{amsmath}
\usepackage{amsfonts}
\usepackage{amssymb}
\usepackage[mathscr]{eucal}
\usepackage{diagrams}
\usepackage{xy,color,graphicx}
\usepackage{hyperref}
\xyoption{all}

\newcommand{\wis}[1]{{\text{\em \usefont{OT1}{cmtt}{m}{n} #1}}}

\newcommand{\N}{\mathbb{N}}

\newcommand{\Af}{\mathbb{A}}
\newcommand{\C}{\mathbb{C}}

\newcommand{\vtx}[1]{*+[o][F-]{\scriptscriptstyle #1}}

\newtheorem{definition}{Definition}

\newtheorem{theorem}{Theorem}

\title{The coordinate biring of $\mathbf{Spec}(\mathbb{Z})/\mathbb{F}_1$}
\author{Lieven Le Bruyn} 
\address{Department of Mathematics, University of Antwerp \\ 
 Middelheimlaan 1, B-2020 Antwerp (Belgium) \\ {\tt lieven.lebruyn@uantwerpen.be}}

\begin{document}
\sloppy

\maketitle

\begin{abstract}
We propose to define $\mathbb{F}_1$-algebras as integral bi-rings with the co-ring structure being the descent data from $\mathbb{Z}$ to $\mathbb{F}_1$. The coordinate bi-ring of $\mathbf{Spec}(\mathbb{Z})/\mathbb{F}_1$ is then the co-ring of integral linear recursive sequences equipped with the Hadamard product. 

We associate a noncommutative moduli space to this setting, show that it is defined over $\mathbb{F}_1$, and has motive $\prod_{n \geq 0} \frac{s-n}{2 \pi}$.
\end{abstract}

\vskip 5mm

Analogous to Borger's approach via $\lambda$-rings \cite{Borger}, we propose in this note $\mathbb{F}_1$-algebras to be integral bi-rings and view the co-ring structure as the descent data from $\mathbb{Z}$ to $\mathbb{F}_1$. This is motivated by the fact that the forgetful functor from integral bi-rings to rings has a right adjoint and hence defines an essential geometric morphism between the corresponding toposes.

Accepting this proposal, we will prove in section~2 that the coordinate bi-ring of $\wis{Spec}(\mathbb{Z})/\mathbb{F}_1$ is  the co-ring of integral linear recursive sequences, equipped with the Hadamard product, as studied by Larson and Taft in \cite{LarsonTaft}.

Over the complex numbers, linear recursive sequences correspond to equivalence classes of canonical linear control systems. Hence one might view them as a noncommutative moduli space, as in \cite{LBReineke}
\[
(M/\mathbb{F}_1)(\C) = \sqcup_n \wis{sys}^c_n . \]
This suggests to look at the larger moduli space of completely controllable and/or completely observable control systems as
\[
(\overline{M}/\mathbb{F}_1)(\C) = \sqcup_n (\wis{sys}^{cc}_n \cup \wis{sys}^{co}_n). \]
We will show in section~3 that this moduli space can be embedded in $\wis{Grass}_2(\infty)$ and has exactly one cell in every possible dimension. Moreover, there is an intriguing involution on this space.

Using a quiver-description of these moduli spaces and Reineke's use of the Harder-Narasinham filtration as in  \cite{Reineke} we will see in section~4 that for every finite field $\mathbb{F}_q$ the number of $\mathbb{F}_q$-points of this noncommutative moduli space is $\sum_{n=0}^{\infty} q^n$, indicating that it is defined over $\mathbb{F}_1$, with corresponding motive
\[
\prod_{n \geq 0} \frac{s-n}{2 \pi}. \]
This observation might be related to Yuri I. Manin's interpretation of Deninger's $\Gamma$-factor at complex arithmetic infinity as the zeta function of (the dual of) infinite dimensional projective space $\mathbb{P}^{\infty}_{\mathbb{F}_1}$, see \cite{Manin1995} and the recent papers \cite[4.3]{Manin2013} and \cite[Intro]{Manin2014}.

\section{Topos theory and $\mathbb{F}_1$-geometry}

Ever since the seminal paper by To\"en and Vaqui\'e \cite{TV}, topos theoretic ideas have been used to develop $\mathbb{F}_1$-geometry, the most recent addition being the 'arithmetic site' discovered by Alain Connes and Katia Consani \cite{CC1}, \cite{CC2}.

Perhaps the most influential of these ideas is due to James Borger \cite{Borger} who proposed to define the category of commutative $\mathbb{F}_1$-algebras to be the category of commutative $\lambda$-rings (without additive torsion). The rationale behind this approach is that Wilkerson \cite{Wilkerson} proved that such rings have a collection of commuting endomorphisms $\psi_p$ which are lifts of the Frobenius-maps. A surprising fact about $\lambda$-rings, see for example \cite{Knutson}, is that the forgetful functor $F$
\[
\xymatrix{\wis{commrings} \ar@/^4ex/[dd]^{W} \ar@/_4ex/[dd]_{U} \\ \\  \lambda-\wis{commrings} \ar[uu]_{F} }
\]
not only has a left adjoint (the free $\lambda$-ring construction $U$), but also a right adjoint, namely the Witt-functor $W$ (see \cite{Hazewinkel} for more details).

Accepting Borger's proposal one then has an essential geometric morphism (see for example \cite[p. 360]{MM}) between the toposes
\[
\xymatrix{\wis{Spaces}_{\mathbb{Z}} \ar@/^4ex/[dd]^{v_{\ast}} \ar@/_4ex/[dd]_{v_{!}} \\ \\  \wis{Spaces}_{\mathbb{F}_1} \ar[uu]_{v^*} }
\]
That is, the base-change functor $v^*$ strips off the $\lambda$-ring structure (which hence can be seen as descent data), the base-forgetting functor $v_!$ sends a space to its Witt-space with natural $\lambda$-structure and the Weil-descent functor $v_*$ sends a space to its arithmetic jet space with natural $\lambda$-structure, see \cite[p.7]{Borger}. In particular, in this proposal $\wis{Spec}(\mathbb{Z})/\mathbb{F}_1$ is the $\mathbb{F}_1$-space corresponding to the $\lambda$-ring of Witt vectors $W(\mathbb{Z})$.

Clearly, one can extend this idea to other categories of (not necessarily commutative) rings for which the forgetful functor has a right adjoint. In this note we consider the case of bi-rings, resp. commutative and cocommutative bi-rings.

With $\wis{rings}$ we denote the category of not necessarily commutative rings {\em having no additive torsion} and ring-morphisms. That is, a ring $A$ has structural multiplication and unit maps
\[
A \otimes_{\mathbb{Z}} A \rTo^m A \qquad \qquad \mathbb{Z} \rTo^u A \]
satisfying the usual commuting diagrams. Dually a co-ring is a $\mathbb{Z}$-module without additive torsion $C$ equipped with co-multiplication and co-unit maps
\[
C \rTo^{\Delta} C \otimes_{\mathbb{Z}} C \qquad \qquad C \rTo^{\epsilon} \mathbb{Z} \]
satisfying the dual commuting diagrams. Between algebras and coalgebras over a field there is the usual Kostant duality, see \cite[Chp. VI]{Sweedler}. We can extend it to the present setting provided we are careful in defining the dual co-ring.

\begin{definition} The dual co-ring $A^o$ of $A \in \wis{rings}$ is the submodule of all $f \in A^*=Hom_{\mathbb{Z}}(A,\mathbb{Z})$ such that $Ker(f)$ contains a twosided ideal $I \triangleleft A$ with $A/I$ a free $\mathbb{Z}$-module of finite rank.
\end{definition}

As $A^o$ is the direct limit of all appropriate $(A/I)^*$, co-multiplication and co-unit on $A^o$ are induced by dualizing the ring structure on $A/I$. With this definition we can copy the proof of \cite[Thm. 6.0.5]{Sweedler}, see \cite{LBmeasuring} for details. That is, for every $A \in \wis{rings}$ and $C \in \wis{co-rings}$ we have a natural one-to-one correspondence
\[
\wis{rings}(A,C^*) \leftrightarrow \wis{co-rings}(C,A^o) \]
Hence, for any torsion free $\mathbb{Z}$-module $V$, the dual of the tensor-algebra on the dual module $T_{\mathbb{Z}}(V^*)^o$ is the cofree co-ring on $V^{**}$. That is, for every co-ring $C$ and additive morphism $f~:~C \rTo V^{**}$ there is a unique co-ring morphism $g~:~C \rTo T(V^*)^o$ making the diagram commute
\[
\xymatrix{C \ar[rr]^g \ar[rrd]_f & & T_{\mathbb{Z}}(V^*)^o \ar[d]^{\pi} \\
& & V^{**}} \]
with $\pi$ dual to the natural inclusion $V^* \rInto T_{\mathbb{Z}}(V^*)$. Indeed, the dual map $V^* \rInto V^{***} \rTo C^*$ gives rise to a ring-morphism $T_{\mathbb{Z}}(V^*) \rTo C^*$ and hence by the above duality to a co-ring map $C \rTo T_{\mathbb{Z}}(V^*)^o$.

Define $H(V)$ to be the union of all sub co-rings $E \rInto T_{\mathbb{Z}}(V^*)^o$ such that $\pi(E) \subset V$, then $H(V)$ is the free co-ring on $V$. In case $V=A \in \wis{rings}$, $H(A)$ becomes a bi-ring with multiplication and unit determined by the maps (induced from the universal property)
\[
\xymatrix{H(A) \otimes_{\mathbb{Z}} H(A) \ar@{.>}[rr] \ar[d]^{\pi \otimes \pi} & & H(A) \ar[d]^{\pi} \\ A \otimes A \ar[rr]_m & & A} \qquad 
\xymatrix{& H(A) \ar[rd]^{\pi} \\ \mathbb{Z} \ar@{.>}[ru] \ar[rr]_i & & A} \]
Similarly, defining $C(V)$ to be the union of all cocommutative co-rings $E \rInto T_{\mathbb{Z}}(V^*)^o$ such that $\pi(E) \subset V$ we have that for $V=A \in \wis{commrings}$ that $C(A)$ is a commutative and cocommutative bi-ring. If $U(-)$ and $V(-)$ are the corresponding free co-ring constructions we have two triples of adjoint functors.

\begin{theorem} The forgetful functors $F$ below have both a left- and a right-adjoint as indicated
\[
\xymatrix{\wis{rings} \ar@/^4ex/[dd]^{H} \ar@/_4ex/[dd]_{U} \\ \\  \wis{bi-rings} \ar[uu]_{F} } \qquad \qquad
\xymatrix{\wis{commrings} \ar@/^4ex/[dd]^{C} \ar@/_4ex/[dd]_{V} \\ \\  \wis{comm.cocomm.bi-rings} \ar[uu]_{F} }
\]
giving rise to essential geometric morphisms
\[
\xymatrix{\wis{NC-Spaces}_{\mathbb{Z}} \ar@/^4ex/[dd]^{v_{\ast}} \ar@/_4ex/[dd]_{v_{!}} \\ \\  \wis{NC-Spaces}_{\mathbb{F}_1} \ar[uu]_{v^*} } \qquad \qquad
\xymatrix{\wis{Spaces}_{\mathbb{Z}} \ar@/^4ex/[dd]^{v_{\ast}} \ar@/_4ex/[dd]_{v_{!}} \\ \\  \wis{Spaces}_{\mathbb{F}_1} \ar[uu]_{v^*} }
\]
\end{theorem}

\section{The coordinate bi-ring of $\wis{Spec}(\mathbb{Z})/\mathbb{F}_1$}

The only commutative ring for which $H(A)=C(A)$ is $A=\mathbb{Z}$. That is, $\wis{Spec}(\mathbb{Z})/\mathbb{F}_1$ will look the same in the commutative and non-commutative world of spaces over $\mathbb{F}_1$. We will now describe its corresponding coordinate bi-ring.

Clearly, $T_{\mathbb{Z}}(\mathbb{Z}) \simeq \mathbb{Z}[t]$ with the natural embedding $\mathbb{Z} \rInto T_{\mathbb{Z}}(\mathbb{Z})$ given by the map $a \mapsto at$. We can identify the full linear dual $\mathbb{Z}[t]^*$ with the module of all infinite sequences $f=(f_n)_{n=0}^{\infty}$ where $f(t^n)=f_n$.

An ideal $I \triangleleft \mathbb{Z}[t]$ such that the quotient $\mathbb{Z}[t]/I$ is without additive torsion is generated by a monic polynomial, $I=(m(t))$ where
\[
m(t) = t^r-a_1 t^{r-1} - \hdots - a_r \]
A sequence $f=(f_n)n \in \mathbb{Z}[t]^*$ such that $I \subset Ker(f)$ has the property that $f(x^k m(t)) = 0$ for all $k \geq 0$, that is,
\[
f_n = a_1 f_{n-1} + a_2 f_{n-2} + \hdots + a_r f_{n-r} \quad \text{for all $n \geq r$} \]
Sequences $f$ satisfying such a relation for some $r >0$ and some $a_i \in \mathbb{Z}$ are called {\em integral linear recursive sequences}. 

Conversely, if $f$ is linear recursive and satisfies the above recursion relation, then $f(\mathbb{Z}[t]m(t))=0$, that is, $f \in \mathbb{Z}[t]^o$. That is, as a module we can identify $H(\mathbb{Z})=C(\mathbb{Z})=\mathbb{Z}[t]^o$ with the module of all integral linear recursive sequences. Observe that the dual map to the inclusion $\mathbb{Z} \rInto \mathbb{Z}[t]$ given by $1 \mapsto t$ gives the projection
\[
\pi~:~\mathbb{Z}[t]^o \rOnto \mathbb{Z} \qquad f=(f_n)_{n=0}^{\infty} \mapsto f_1 \]
The co-ring structure on $\mathbb{Z}[t]^o$ is determined by the (usual) multiplication on $\mathbb{Z}[t]$. The coalgebra structure of $k[t]^o$ was described in \cite{PetersonTaft} and \cite{LarsonTaft} in case $k$ is an (algebraically closed) field. We will recall the result for $\mathbb{Q}[t]^o$.

Start with a rational linear recursive sequence $(f_n)_n$ and let $t$ be the largest integer such that the rows of the symmetric $t \times t$ matrix
\[
H(f) = \begin{bmatrix} f_0 & f_1 & f_2 & \hdots & f_{t-1} \\
f_1 & f_2 & f_3 & \hdots & f_t \\
f_2 & f_3 & f_4 & \hdots & f_{y+1} \\
\vdots & \vdots & \vdots & \cdots & \vdots \\
f_{t-1} & f_t & f_{t+1} & \hdots & f_{2t-2} \end{bmatrix} \]
are linearly independent. The matrix $H(f)$ is called the {\em Hankel matrix} of $f$ and is by definition invertible with inverse $H^{-1} = (s_{ij})_{0 \leq i,j \leq t-1}$. If we define $D^if$ to be the rational recursive sequence $(D^if)_n = f_{n+i}$, then the coproduct on $\mathbb{Q}[t]^o$ is given by
\[
\Delta(f) = \sum_{i,j=0}^{t-1} s_{ij} (D^if) \otimes (D^jf) \]
Clearly, if for $f \in \mathbb{Z}[t]^o$ the determinant of the Hankel matrix $H(f)$ is $\pm 1$, the same formula applies to $\Delta(f)$. But, in general $\Delta(f)$ cannot be diagonalized in terms of $f, Df, \hdots$ with integer coefficients and we have no other option but to describe the comultiplication on $\mathbb{Z}[t]^o$ as the inductive limit of the coproducts determined by taking the full linear duals of the form $(\mathbb{Z}[t]/(m(t)))^*$ where $m(t)$ runs over the monic integral polynomials and the limit is taken with respect to divisibility of monic polynomials.

As a consequence we also have to describe the dual ring $(\mathbb{Z}[t]^o)^*$ as the inverse limit
\[
(\mathbb{Z}[t]^o)^* = \underset{\leftarrow}{{\tt lim}}~\frac{\mathbb{Z}[t]}{(m(t))} \]
with respect to the divisibility relation on monic polynomials. This completion at all monic polynomials has been studied by Habiro, \cite{Habiro}. The subring which is the completion at all cyclic polynomials has appeared before in $\mathbb{F}_1$-geometry, see for example \cite{ManinHabiro}.

\begin{theorem} The coordinate bi-ring of $\wis{Spec}(\mathbb{Z})/\mathbb{F}_1$ is the co-ring of all integral linear recursive sequences, equipped with the Hadamard product.
\end{theorem}

\begin{proof}
As multiplication and unit on $\mathbb{Z}[t]^o$ are uniquely determined by the commuting diagrams
\[
\xymatrix{\mathbb{Z}[t]^o \otimes_{\mathbb{Z}} \mathbb{Z}[t]^o \ar@{.>}[rr] \ar[d]^{\pi \otimes \pi} & & \mathbb{Z}[t]^o \ar[d]^{\pi} \\ \mathbb{Z} \otimes \mathbb{Z} \ar[rr]_m & & \mathbb{Z}} \qquad 
\xymatrix{& \mathbb{Z}[t]^o \ar[rd]^{\pi} \\ \mathbb{Z} \ar@{.>}[ru] \ar[rr]_i & & \mathbb{Z}} \]
and as $\pi(f)=f_1$ we deduce that the product on $\mathbb{Z}[t]^o$ is the {\em Hadamard product}, see \cite{LarsonTaft}, given by $f.g = (f_i.g_i)_i$ and that the unit is the sequence $(1,1,1,\hdots)$.
\end{proof}

The bi-ring $\mathbb{Z}[t]^o$ has a lot of additional structure. For example, the family of commuting ring-morphisms $\{ \psi^n~:~n \in \N_+ \}$ determining the $\lambda$-structure on $\mathbb{Z}[t]$ (that is, $\psi^n(t) = t^n$) induce via the natural inclusion
\[
\xymatrix{\mathbb{Z}[t] \ar[rr]^{\psi^n} \ar@{.>}[rrd]_{\tilde{\psi}^n} & & \mathbb{Z}[t] \ar[d] \\
& & \mathbb{Z}[t]^{o \ast} }
\]
a commuting family of ring maps $\tilde{\psi}^n$. By Kostant duality
\[
\wis{rings}(\mathbb{Z}[t],\mathbb{Z}[t]^{o \ast}) \simeq \wis{co-rings}(\mathbb{Z}[t]^o,\mathbb{Z}[t]^o) \]
this gives a commuting family of co-ring maps on $\mathbb{Z}[t]^o$, that is, an action of $\mathbb{N}_+^{\times}$ on the co-ring $\mathbb{Z}[t]^o$. Further, by the adjointness
\[
\wis{bi-rings}(\mathbb{Z}[t]^o,\mathbb{Z}[t]^o) \simeq \wis{rings}(\mathbb{Z}[t]^o,\mathbb{Z}) \]
the family of projection ring-maps $\pi_n : \mathbb{Z}[t]^o \rOnto \mathbb{Z}$ sending $\pi_n((f_m)m) \mapsto f_n$, gives a family of bi-ring endomorphisms on $\mathbb{Z}[t]^o$.

\section{A noncommutative moduli space}

In several proposals to $\mathbb{F}_1$-geometry, noncommutative geometric objects seem to appear quite naturally. In Borger's proposal \cite{Borger}, one quickly encounters the Bost-Connes system \cite{CC3}. 

In the present proposal, the geometric object associated to $\mathbb{Z}[t]^o$ will be a noncommutative moduli space, in the sense of Maxim Kontsevich \cite{KontFormal}, which consists of gluing together ordinary moduli spaces into an infinite dimensional variety, which is then controlled by a noncommutative formally smooth algebra, see \cite{LBReineke} for more details.

Recall from \cite[Part VI-VII]{Tannenbaum} that a linear recursive sequence with Hankel matrix of size $n$ can be realized as the input-output behavior of a {\em canonical} (that is, completely controllable and completely observable) linear dynamical system with one input, one output and an $n$ dimensional state space. In \cite{LBReineke} the corresponding moduli spaces (with varying $n$) were shown to form a noncommutative moduli space in the above sense, even if we allow arbitrary but fixed input- and output-dimensions. 

We will quickly run through the definitions and the proof of \cite{LBReineke} in the special case of interest here, identifying this noncommutative manifold with a specific cell-subcomplex of $\wis{Grass}_2(\infty)$.

A {\em linear control system} $\Sigma$ with one input and one output and $n$-dimensional state space is determined by the system of linear differential equations
\[
\begin{cases}
\frac{d x}{dt} &= Ax + Bu \\ y &= Cx
\end{cases} \]
where $u(t) \in \C$ is the input (or {\em control}) at time $t$, $x(t) \in \C^n$ is the {\em state} of the system and $y(t) \in \C$ is its {\em output}.

That is, $\Sigma$ is determined by a triple $(A,B,C)$ where $A \in M_n(\C)$, $B \in \C^n$ and $C \in \C^{n \ast}$. Another system $\Sigma'=(A',B',C')$ is said to be equivalent to $\Sigma$ if there is a base-change matrix $g \in GL_n(\C)$ such that
\[
A' = g A g^{-1} \qquad B' = g B \quad \text{and} \quad C' = C g^{-1} \]
$\Sigma$ is said to be {\em completely controllable} (resp. {\em completely observable}) if and only if the $n \times n$ matrix
\[
c(\Sigma) = \begin{bmatrix} B & AB & A^2B & \hdots & A^{n-1} B \end{bmatrix} \qquad (\text{resp.} \quad o(\Sigma) = \begin{bmatrix} C \\ CA \\ CA^2 \\ \vdots \\ CA^{n-1} \end{bmatrix} \quad) \]
is invertible. These conditions define $GL_n$-open subsets $V^{cc}_n$ (resp. $V^{co}_n$) of $V_n = M_n(\C) \oplus \C^n \oplus \C^{n \ast}$ and the corresponding orbit spaces
\[
\wis{sys}^{cc}_n = V^{cc}_n / GL_n \quad \text{and} \quad \wis{sys}^{co}_n = V^{co}_n / GL_n \]
are known to be smooth quasi-projective varieties of dimension $2n$, see for example \cite[Part IV]{Tannenbaum}, and isomorphic to each other sending a system $\Sigma=(A,B,C)$ to the transposed system $\Sigma^t = (A^t,C^t,B^t)$.

$\Sigma$ is said to be {\em canonical} if it is both completely controllable and completely observable, and the corresponding moduli space
\[
\wis{sys}^c_n = (V^{cc}_n \cap V^{co}_n) / GL_n \]
classifies canonical systems such that the linear recursive sequences $(f_i = CA^iB)_{i \in \N}$ are equal, see \cite[Part VI-VII]{Tannenbaum}. Conversely, a linear recursive sequence $(f_i)_{i \in \N}$ with Hankel matrix of size $n$ is realizable, that is $f_i = C A^i B$ for all $i \in \N$, by a canonical system $\Sigma$ in $\wis{sys}_n^c$. 

In the next section we will show that these moduli spaces can be defined over $\mathbb{F}_1$. For this reason we define a noncommutative moduli space $M/\mathbb{F}_1$ and its 'closure' $\overline{M}/\mathbb{F}_1$ having as its $\C$-points
\[
(M / \mathbb{F}_1)(\C) = \bigsqcup_n \wis{sys}^c_n  \quad \text{and} \quad (\overline{M}/ \mathbb{F}_1)(\C) = \bigsqcup_n (\wis{sys}^{cc}_n \cup \wis{sys}^{co}_n).  \]
As $\mathbb{Z}/\mathbb{F}_1 = \mathbb{Z}[t]^o$, it only takes a leap of faith to envision that $(\overline{M}/\mathbb{F}_1)(\mathbb{C})$ might be related to $\overline{\wis{Spec}(\mathbb{Z})}/\mathbb{F}_1$ at complex arithmetic infinity.

We will now embed both of these noncommutative moduli spaces into $\wis{Grass}_2(\infty)$.
A completely controllable system $\Sigma = (A,B,C)$ is equivalent to one in canonical form, and if we assign to a system the $n \times n+2$ matrix $M_{\Sigma}=(B,C^{\tau},A)$, then the matrix corresponding to a system in canonical form is
\[
M_{\Sigma} = \begin{bmatrix} 1 & c_1 & 0 & \hdots & 0 & a_1 \\
0 & c_2 & 1 & & 0 & a_2 \\
\vdots & \vdots & & \ddots & & \vdots \\
0 & c_n & 0 &  & 1 & a_n \end{bmatrix} \]
which determines a point in $\wis{Grass}_n(n+2)$. Under the duality with $\wis{Grass}_2(n+2)$ this matrix corresponds to the $2 \times n$ matrix
\[
K_{\Sigma} = \begin{bmatrix} -c_1 & 1 & -c_2 & \hdots & -c_n & 0 \\ -a_1 & 0 & -a_2 & \hdots & -a_n & 1 \end{bmatrix} \]
which is a point in the standard open cell of dimension $2n$ with multi-index $\{ 2, n+2 \}$. Similarly, any equivalence class of completely observable system corresponds to a point in the standard cell of dimension $2n$ in $\wis{Grass}_2(n+2)$ with multi-index $\{ 1, n+2 \}$. The complement of the intersection of these two cells in either of them  is a $2n-1$-dimensional cell.

\begin{theorem} For every $n$ we have a cell decomposition of the moduli space $\wis{sys}^{cc}_n \cup \wis{sys}^{co}_n$ as $\C^{2n} \sqcup \C^{2n-1}$ in the Grassmannian $\wis{Grass}_2(n+2)$ as
\[
\begin{bmatrix} \ast & 1 & \ast & \hdots & \ast & 0 \\ \ast & 0 & \ast & \hdots & \ast & 1 \end{bmatrix} \sqcup
\begin{bmatrix} 1 & 0 & \ast & \hdots & \ast & 0 \\ 0 & \ast & \ast & \hdots & \ast & 1 \end{bmatrix} \]
As a consequence, we have a cell decomposition of $(\overline{M}/\mathbb{F}_1)(\C) = \sqcup_n (\wis{sys}^{cc}_n \cup \wis{sys}^{co}_n)$ in $\wis{Grass}_2(\infty)$ as
\[
\begin{bmatrix} \ast & 1 & \ast & \hdots & \ast & 0 & 0 & 0 & \hdots \\
\ast & 0 & \ast & \hdots & \ast & 1 & 0 & 0 & \hdots \end{bmatrix} \sqcup
\begin{bmatrix} 1 & 0 & \ast & \hdots & \ast & 0 & 0 & 0 & \hdots \\
0 & \ast & \ast & \hdots & \ast & 1 & 0 & 0 & \hdots \end{bmatrix} \]
that is, all matrices of this form, with a tail of zero columns, the last non-zero column $x$ determining the dimension $n=x-2$ of the state space of the corresponding  system.
\end{theorem}

Observe that sending a system $\Sigma=(A,B,C)$ to its transposed system $\Sigma^t=(A^t,C^t,B^t)$ induces an involution on this space.

\section{Manin's motivic mystery}

In \cite{Manin1995}, Yuri I. Manin suggested the existence of a category of $\mathbb{F}_1$-motives visible through the $q=1$ point count of $\mathbb{F}_1$-schemes, and, this prediction was later justified by Chr. Soul\'e's version of $\mathbb{F}_1$-geometry, \cite{Soule}.

In his lecture notes on zeta functions and motives \cite{Manin1995}, Yuri I. Manin defined the zeta function of $\mathbb{P}^k_{\mathbb{F}_1}$, the $k$-dimensional projective space over the hypothetical field with one element $\mathbb{F}_1$, to be
\[
\prod_{n=0}^k \frac{s-n}{2 \pi}. \]
In \cite{Deninger1991}, Christopher Deninger represented the basic $\Gamma$-factor at complex arithmetic infinity as the infinite determinant of the Frobenius map and a regularized product
\[
\Gamma_{\mathbb{C}}(s)^{-1} = \frac{(2 \pi)^s}{\Gamma(s)} = \prod_{n \geq 0} \frac{s+n}{2 \pi}. \]
Comparing both formulas led Manin to suggest that this $\Gamma$-factor might be viewed as the zeta-function (or motive) of the (dual of) infinite dimensional projective space over $\mathbb{F}_1$, see also his recent papers \cite[4.3]{Manin2013} and \cite[Intro]{Manin2014}.
As such, one might expect this local factor to appear naturally in the zeta-function of $\overline{\wis{Spec}(\mathbb{Z})}/\mathbb{F}_1$.

There is a close analogy between Manin's $\mathbb{F}_1$-motives and Kurokawa's zeta functions of schemes over $\mathbb{F}_1$, \cite{Kurokawa}. N. Kurokawa defined an algebraic set $X$ defined over $\mathbb{Z}$ to be of $\mathbb{F}_1$-type if there exists a  polynomial with integer coefficients
\[
N_X(t) = \sum_{k=0}^n a_k t^k \qquad \text{such that} \qquad \# X(\mathbb{F}_q) = N_X(q) \]
for all finite fields $\mathbb{F}_q$. In this case, Kurokawa defines the zeta-function of $X$ over $\mathbb{F}_1$ to be the expression
\[
\zeta(s,X/\mathbb{F}_1) = \prod_{k=0}^n \frac{1}{(s-k)^{a_k}} \]
see also Soul\'e \cite{Soule} and Deitmar \cite{Deitmar} for closely related definitions. For example, for $\Af^n/\mathbb{F}_1$ (resp. the $n$-th power of the Tate motive $\mathbf{T}$) and $\mathbb{P}^n/\mathbb{F}_1$ we have the following expressions of Kurokawa's zeta function versus Manin's $\mathbb{F}_1$-motive
\[
\begin{cases}
\zeta(s,\mathbb{A}^n/\mathbb{F}_1) = \frac{1}{s-n} & Z(\mathbf{T}^{\times n},s) = \frac{s-n}{2 \pi} \\ \\ \zeta(s,\mathbb{P}^n/\mathbb{F}_1) = \prod_{k=0}^n \frac{1}{s-k} & Z(\mathbb{P}^n_{\mathbb{F}_1},s) = \prod_{k=0}^n \frac{s-k}{2 \pi}
\end{cases}
 \]
We will now show that each of the components $\wis{sys}_n^{cc} \cup \wis{sys}^{cc}_n$ is of $\mathbb{F}_1$-type and compute its zeta-function and $\mathbb{F}_1$-motive.

To a linear control system $\Sigma = (A,B,C)$ with $n$-dimensional state space (and one-dimensional input and output) we associate the representation of the quiver $Q$
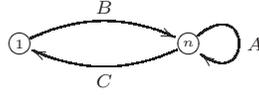
\[
\xymatrix{ \vtx{1} \ar@/^2ex/[rr]^B & & \vtx{n} \ar@/^2ex/[ll]^C \ar@(ur,dr)^A} \]
of dimension vector $\alpha_n = (1,n)$. If we take the stability structures $\theta_+ = (-n,1)$ and $\theta_-=(n,-1)$, then it follows from \cite[Lemma 1 \& 2]{LBReineke} that
\begin{itemize}
\item{Equivalence classes of canonical systems correspond to isomorphism classes of simple quiver-representations.}
\item{Equivalence classes of completely controllable  systems correspond to isomorphism classes of $\theta_+$-stable quiver-representations.}
\item{Equivalence classes of completely observable systems correspond to isomorphism classes of $\theta_-$-stable quiver-representations.}
\end{itemize}
That is, if $\wis{simple}_{\alpha_n}(Q)$ denotes the quasi-affine quotient variety of isomorphism classes of simple $\alpha_n$-dimensional $Q$-representations, and if $\wis{moduli}^{\theta_\pm}_{\alpha_n}(Q)$ denote the quiver moduli spaces of $\theta_{\pm}$-semistable representations, as introduced and studied by A. King \cite{King}, then we have isomorphisms of varieties
\[
\wis{sys}^c_n \simeq \wis{simple}_{\alpha_n}(Q), \qquad \wis{sys}^{cc}_n \simeq \wis{moduli}^{\theta_+}_{\alpha_n}(Q), \qquad \wis{sys}^{co}_n \simeq \wis{moduli}^{\theta_-}_{\alpha_n}(Q). \]

In \cite{Reineke} the Harder-Narasinham filtration was used to compute the cohomology of such quiver moduli spaces (at least when the quiver has no oriented cycles) and the same methods apply to compute the number of $\mathbb{F}_q$-points of quiver moduli spaces in general. From \cite{LBReineke} we obtain
\[
\begin{cases}
\#~\wis{sys}^{co}_n(\mathbb{F}_q) = \#~\wis{moduli}^{\theta_+}_{\alpha_n}(Q)(\mathbb{F}_q) = q^{2n} \\ \\ 
\#~\wis{sys}^{cc}_n(\mathbb{F}_q) = \#~\wis{moduli}^{\theta_-}_{\alpha_n}(Q)(\mathbb{F}_q) = q^{2n} \\ \\
\#~\wis{sys}^c_n(\mathbb{F}_q) = \#~\wis{simple}_{\alpha_n}(Q)(\mathbb{F}_q) = q^{2n} -q ^{2n-1}
\end{cases}
\]
and therefore $\wis{sys}^{co}_n \cup \wis{sys}^{cc}_n$ is of $\mathbb{F}_1$-type and we have
\[
\#~(\wis{sys}^{co}_n \cup \wis{sys}^{cc}_n)(\mathbb{F}_q) = q^{2n} + q^{2n-1} \]
giving us their Kurokawa zeta-function and Manin motive over $\mathbb{F}_1$.

\begin{theorem} $(\overline{M}/\mathbb{F}_1)(\C)$ is of $\mathbb{F}_1$-type with corresponding $\mathbb{F}_1$-motive
\[
\prod_{k=0}^{\infty} \frac{s-k}{2 \pi} \]
\end{theorem}


\begin{thebibliography}{10}

\bibitem{Borger}
James Borger, {\it $\Lambda$-rings and the field with one element}, {\tt arXiv:0906.3146} (2009)


\bibitem{CC3}
Alain Connes and Katia Consani, {\it On the arithmetic of the BC-system}, {\tt arXiv:1103.4672} (2011)

\bibitem{CC1}
Alain Connes and Katia Consani, {\it The arithmetic site}, {\tt arXiv:1405.4527} (2014)

\bibitem{CC2}
Alain Connes and Katia Consani, {\it Geometry of the arithmetic site}, {\tt arXiv:1502.05580} (2015)

\bibitem{Deitmar}
A. Deitmar, {\it Schemes over $\mathbb{F}_1$}, {\tt arXiv:0404185} (2004)

\bibitem{Deninger1991}
Christopher Deninger, {\it On the $\Gamma$-factors attached to motives}, Invent. Math. 104 (1991) 245-261

\bibitem{Habiro}
Kazuo Habiro, {\it Cyclotomic completions of polynomials}, {\tt arXiv:0209.324} (2002)

\bibitem{Hazewinkel}
Michiel Hazewinkel, {\it Witt vectors, part 1}, {\tt arXiv:0804.3888} (2008)

\bibitem{King}
Alastair King, {\it Moduli of representations of finite dimensional algebras}, Quat. J. Math. Oxford Ser. (2) {\bf 45} (1994) 515-530

\bibitem{Knutson}
Donald Knutson, {\it $\lambda$-rings and the representation theory of the symmetric group}, Springer Lect. Notes in Math. 308 (1973)

\bibitem{KontFormal}
Maxim Kontsevich, {\it Formal non-commutative symplectic geometry}, Gelfand Seminar 1990-1992, Birkh\"auser (1993) 173-187

\bibitem{Kurokawa}
Nobushige Kurokawa, {\it Zeta functions over $\mathbb{F}_1$}, Proc. Japan Acad., 81, Ser. A (2005) 180-184

\bibitem{LarsonTaft}
Richard G. Larson and Earl J. Taft, {\it The algebraic structure of linearly recursive sequences under Hadamard product}, Israel Journal of Mathematics,
Volume 72 (1990) 118-132

\bibitem{LBmeasuring}
Lieven Le Bruyn, {\it Universal bialgebras associated with orders},  Comm. Algebra, 10 (1982) 457-478


\bibitem{LBReineke}
Lieven Le Bruyn and Markus Reineke, {\it Canonical systems and non-commutative geometry}, {\tt arXiv:math.AG/0303304} (2003)

\bibitem{MM}
Saunders Mac Lane and Ieke Moerdijk, {\it Sheaves in geometry and logic, a first introduction to topos theory}, Universitext Springer-Verlag (1992)

\bibitem{Manin1995}
Yuri I. Manin, {\it Lectures on zeta functions and motives (according to Deninger and Kurokawa). In: Columbia University Number Theory Seminar (1992)}, Ast\'erisque 228 (1995) 121-164

\bibitem{ManinHabiro}
Yuri I. Manin, {\it Cyclotomy and analytic geometry over $\mathbb{F}_1$}. In: Quanta of Maths. Conference in honour of Alain Connes. Clay Math. Proceedings, vol 11 (2010), 385-408, {\tt arXiv:math.AG/0809.2716}

\bibitem{Manin2013}
Yuri I. Manin, {\it Numbers as functions}, {\tt arXiv:1312.5160} (2013)

\bibitem{Manin2014}
Yuri I. Manin, {\it Local zeta functions and geometries under $\wis{Spec}(\mathbb{Z})$}, {\tt arXiv:1407.4969} (2014)

\bibitem{PetersonTaft}
Brian Peterson and Earl J. Taft, {\it The Hopf algebra of recursive sequences}, Aequationes Mathematicae {\bf 20} (1980) 1-17

\bibitem{Reineke}
Markus Reineke, {\it The Harder-Narasinham system in quantum groups and cohomology of quiver moduli}, Invent. Math. {\bf 152} (2003) 349-368

\bibitem{Soule}
Christophe Soul\'e, {\it Les vari\'et\'es sur le corps \`a un \'el\'ement}, Mosc. Math. J. {\bf 4} (2004) 217-244 

\bibitem{Sweedler}
Moss E. Sweedler, {\it Hopf Algebras}, W.A. Benjamin Inc., New York (1969)


\bibitem{Tannenbaum}
Allen Tannenbaum, {\it Invariance and system theory: algebraic and geometric aspects}, Lect. Notes in Math. 845 (1981) Springer-Verlag

\bibitem{TV}
B. To\"en and M. Vaqui\'e, {\it Au-dessous de $\wis{Spec}(\mathbb{Z})$}, Journal of K-theory 3 (2009) 437-500, {\tt arXiv:0509684}

\bibitem{Wilkerson}
Clarence Wilkerson, {\it Lambda-rings, binomial domains and vector bundles over $CP(\infty)$}, Comm. Algebra 10(3) (1982) 311-328

\end{thebibliography}
\end{document}